\documentclass[12pt,leqno]{article}
\usepackage{amssymb,amsmath,amsthm}

\def\beq{\begin{equation}}
\def\eeq{\end{equation}}

\def\ha{\hat{a}}
\def\hb{\hat{b}}
\def\hR{\hat{R}}

\def\bC{\mathbf{C}}
\def\bI{\mathbf{I}}

\def\be{\mathbf{e}}
\def\bn{\mathbf{n}}

\def\cF{\mathcal{F}}

\def\eps{\varepsilon}

\newtheorem{theorem}{Theorem}
\newtheorem{lemma}[theorem]{Lemma}

\begin{document}

\title{Fitting circles to scattered data: parameter estimates have no moments}

\author{ N. Chernov\\
Department of Mathematics\\ University of Alabama at Birmingham\\
Birmingham, AL 35294\\ chernov@math.uab.edu\\ Fax 1-205-934-9025}

\date{}

\maketitle

\begin{abstract}
We study a nonlinear regression problem of fitting a circle (or a
circular arc) to scattered data. We prove that under any standard
assumptions on the statistical distribution of errors that are
commonly adopted in the literature, the estimates of the circle
center and radius have infinite moments. We also discuss
methodological implications of this fact.
\end{abstract}

\begin{center}
\end{center}

\section{Introduction}
\label{secI}

Regression models in which all variables are subject to errors are
known as error-in-variables (EIV) models. The EIV regression problem
is quite different (and far more difficult) than the classical
regression where the independent variable is assumed to be
error-free. The EIV regression, even in the linear case, presents
extremely challenging questions and leads to some counterintuitive
results (some of them are mentioned below).

This work is devoted to a nonlinear EIV model where one fits a
circle to scattered data. This is one of the basic tasks in pattern
recognition and computer vision. The need of fitting circles to
planar images also arises in biology and medicine, nuclear physics,
archeology, industry, and other areas of human practice.

The most popular method used to solve this problem is orthogonal
least squares, i.e.\ the minimization of the sum of squares of the
distances from the data points to the fitting contour. This method
is often called \emph{geometric fit} or \emph{orthogonal distance
regression} (ODR).

Fitting a circle to observed points $(x_1, y_1), \ldots, (x_n, y_n)$
amounts to minimizing the objective function
\beq
       \cF(a,b,R) =  \sum_{i=1}^n \bigl[
       \sqrt{(x_i-a)^2+(y_i-b)^2} - R\bigr]^2,
         \label{Fmain}
\eeq
where $(a,b)$ denotes the center and $R$ the radius of the circle.
Then the parameters of the best fitting circle are defined by
\beq  \label{argmin}
   (\ha, \hb, \hR) = \, {\rm argmin}\, \cF(a,b,R).
\eeq
To explore the statistical properties of the estimates $\ha, \hb,
\hR$ one needs to make assumptions on the probability distribution
of the data points. It is commonly assumed that each $(x_i, y_i)$ is
a noisy observation of some \emph{true point} $(x^{\ast}_i,
y^{\ast}_i)$, i.e.
\beq \label{errors}
   x_i = x^{\ast}_i + \delta_i, \qquad
   y_i = y^{\ast}_i + \eps_i,
   \qquad i=1,\ldots,n,
\eeq
where $(\delta_1, \eps_1), \ldots, (\delta_n, \eps_n)$ are $n$
independent random vectors, usually with zero mean.

A standard assumption is that each $(\delta_i, \eps_i)$ is a normal
(Gaussian) vector with some covariance matrix $\bC_i$. The simplest
choice is $\bC_i = \sigma^2 \bI$, in which case all errors
$\eps_i$'s and $\delta_i$'s are i.i.d.\ normal random variables with
zero mean and a common variance $\sigma^2$. In that case the
geometric fit \eqref{argmin} coincides with the maximum likelihood
estimate (MLE), see Chan 1965.

The true points $(x^{\ast}_i, y^{\ast}_i)$ are supposed to lie on a
`true circle', i.e.\ satisfy
\beq
  (x^{\ast}_i - a^{\ast})^2 + (y^{\ast}_i - b^{\ast})^2 = (R^{\ast})^2,
   \qquad i=1,\ldots,n,
\eeq
where $(a^{\ast}, b^{\ast}, R^{\ast})$ denote the `true' (unknown)
parameters. Therefore
$$
   x^{\ast}_i = a^{\ast} + R^{\ast} \cos\varphi_i,
   \qquad y^{\ast}_i = b^{\ast} + R^{\ast} \sin\varphi_i,
$$
where $\varphi_1, \ldots, \varphi_n$ specify the location of the
true points on the true circle.

The angles $\varphi_1, \ldots, \varphi_n$ can be regarded as fixed
unknowns, then they have to be treated as additional parameters of
the model (often called \emph{incidental} or \emph{latent}
parameters). This setup is known as a \emph{functional model}, see
Chan 1965.

Alternatively, $\varphi_1, \ldots, \varphi_n$ can be regarded as
independent realizations of a random variable with a certain
probability distribution on $[0,2 \pi]$; then one gets the so called
\emph{structural model}, see Anderson 1981 or Berman and Culpin
1986. Both models are widely used in the literature.

Many authors study the distribution of the estimates $\ha, \hb, \hR$
under the above assumptions and try to evaluate their biases and
covariance matrix. Our main result is

\begin{theorem}
If the probability distribution of each vector $(\delta_i, \eps_i)$
has a continuous strictly positive density, then $\ha, \hb, \hR$ do
not have moments, i.e.\
$$
   E(|\ha|) = E(|\hb|) = E(\hR) = \infty.
$$
Thus the estimates $\ha, \hb, \hR$ have no mean values or variances.
\end{theorem}

Our assumptions include (but are not limited to) normally
distributed errors. The distribution of $(\delta_i, \eps_i)$ need
not be the same for different $i$'s, it may depend on $i$, but the
vectors $(\delta_i, \eps_i)$ must be independent. The mean value of
$(\delta_i, \eps_i)$ need not be zero. The theorem is valid for
every $n \geq 3$.

\section{Historical remarks} \label{SecHR}

Our result is not entirely surprising as a similar theorem has been
proven for orthogonal least squares lines by Anderson 1976. Suppose
one fits a line $y=\alpha+\beta x$ to data points $(x_1, y_1),
\ldots, (x_n, y_n)$ by minimizing the sum of squares of (orthogonal)
distances, i.e.\ the estimates are defined by
$$
   (\hat{\alpha},\hat{\beta}) = \,{\rm argmin}\,
   \frac{1}{1+\beta^2}\,\sum_{i=1}^n (y_i - \alpha - \beta x_i)^2.
$$
Again the observed points are random perturbations of some true
points, in the sense of \eqref{errors}, which lie on an unknown true
line, i.e.\ satisfy
\beq
  y^{\ast}_i = \alpha^{\ast} + \beta^{\ast}x^{\ast}_i,
   \qquad i=1,\ldots,n.
\eeq
The true points are either fixed parameters (making it a functional
model), or randomly sampled on the true line (structural model).

\begin{theorem}[Anderson 1976]
If the errors $\delta_i$'s and $\eps_i$'s are i.i.d.\ normal random
variables with zero mean and a common variance $\sigma^2>0$, then
$\hat{\alpha}$ and $\hat{\beta}$ do not have moments, i.e.\ $
   E(|\hat{\alpha}|) = E(|\hat{\beta}|) = \infty.
$
\end{theorem}

Until Anderson's discovery, statisticians were used to employ Taylor
expansion to derive some `approximate' formulas for the moments of
the estimates $\hat{\alpha}$ and $\hat{\beta}$ (including their
means and variances). Anderson demonstrated that all those formulas
should be regarded as \emph{moments of some approximations}, rather
than `approximate moments'.

Anderson's result was rather sensational at the time, it was
followed by heated discussions and a period of acute interest in the
linear EIV regression. It also created methodological problems which
we discuss in the next section.

Anderson proved his theorem by using an explicit formula for the
density function of $\hat{\beta}$ (that formula was mentioned but
not given in his paper; it appeared in a later paper by Anderson and
Sawa 1982. Anderson also remarked that his result can be
`intuitively seen' from a well known formula for $\hat{\beta}$:
\beq \label{bline}
   \hat{\beta} = \frac{s_{yy}-s_{xx}+\sqrt{(s_{yy}-s_{xx})^2
   +4s_{xy}^2}}{2s_{xy}},
\eeq
where standard statistical notation are used: $ s_{xx} =
\sum_{i=1}^n (x_i - \bar{x})^2$, $s_{yy} = \sum_{i=1}^n (y_i -
\bar{y})^2$, $s_{xy} = \sum_{i=1}^n (x_i - \bar{x})(y_i - \bar{y})$,
and $\bar{x} = \frac 1n \sum_{i=1}^n x_i$, $\bar{y} = \frac 1n
\sum_{i=1}^n y_i$. Anderson 1976 says that the (continuous) density
of $s_{xx}$, $s_{yy}$, and $s_{xy}$ for which the numerator in
\eqref{bline} is different from 0 and the denominator is equal to 0
is positive, hence the integral of the product of $\hb$ and its
density diverges.

This `intuitive' explanation can be easily converted into a rigorous
proof, and then one readily extends Anderson's theorem to arbitrary
distributions of errors as long as they have continuous strictly
positive densities, like in our Theorem~1. (Alternatively, one can
easily modify our constructions below to achieve this goal; this is
all fairly straightforward, so we omit details.)

We note that Anderson's result was recently extended to some other
estimates of the linear parameters $\alpha$ and $\beta$, see Chen
and Kukush 2006 and an example in Zelniker and Clarkson 2006.

The problem of fitting circles (as well as other nonlinear curves)
to data is technically much more difficult than that of fitting
lines. In particular, there are no explicit formulas for the
estimates $\ha$, $\hb$ or $\hR$, analogous to \eqref{bline}, let
alone explicit formulas for their probability densities. All the
methods of computing the estimates $\ha$, $\hb$ or $\hR$ are based
on iterative numerical schemes.

All this makes the problem of fitting circles (ellipses, etc.) so
much different from that of fitting lines that Anderson's result
apparently passed unnoticed by the `curve fitting community'. It is
still commonly believed that the curve's parameters have finite
moments; thus many researchers try to minimize their bias and
variances or compute Cramer-Rao lower bounds on the covariance
matrix, see e.g.\ Kanatani 1998 or Chernov and Lesort 2004.

Our theorem shows that the true moments are infinite in the case of
fitting circles. We believe this result holds for ellipses and other
types of curves, and we plan to investigate this issue.

We note that Kukush et al.\ 2004 recently modified the geometric
fitting of ellipses to data in order to ensure the consistency of
the parameter estimates (as the orthogonal regression estimates are
inconsistent). They noted that their modified estimates had infinite
moments, which at the time seemed to be a price to pay for
consistency. It is now clear that the lack of moments is a rather
general property of estimates under the EIV model.

\section{Methodological issues} \label{SecFD}

When Anderson proved his Theorem 2, it immediately lead to
fundamental methodological questions: can one trust a statistical
estimate that has an infinite mean square error (not to mention
infinite bias)? Can such an estimate be better than others which
have finite moments?

\medskip \textbf{Fitting lines.}
To explore this issue, Anderson 1976 and 1984, Kunitomo 1980, and
others compared the MLE estimate $\hat{\beta}$ given by
\eqref{bline} with the classical estimate $\hat{\beta} = s_{xy} /
s_{xx}$ of the slope of the regression line that is known to be
optimal when $x_i$'s are error-free (i.e., $\delta_i=0$). They
denote the former by $\hat{\beta}_{\rm M}$ (Maximum likelihood) and
the latter by $\hat{\beta}_{\rm L}$ (Least squares); of course, both
estimates were studied in the framework of the EIV model described
in Section~\ref{SecHR}. Their results can be summarized in two
seemingly conflicting verdicts:
\begin{itemize} \item[(a)] The mean square error of
$\hat{\beta}_{\rm M}$ is infinite, and that of $\hat{\beta}_{\rm L}$
is finite (whenever $n \geq 4$), thus $\hat{\beta}_{\rm L}$ appears
(infinitely!) more accurate;
\item[(b)] The estimate $\hat{\beta}_{\rm M}$ is consistent and
asymptotically unbiased, while $\hat{\beta}_{\rm L}$ is inconsistent
and asymptotically biased (unless $\beta=0$).
\end{itemize}
Besides, Anderson 1976 shows that if $\beta \neq 0$, then
$$
   P\bigl(|\hat{\beta}_{\rm M}-\beta|>t\bigr) <
   P\bigl(|\hat{\beta}_{\rm L}-\beta|>t\bigr)
$$
for all $t>0$ of practical interest, i.e.\ the accuracy of
$\hat{\beta}_{\rm M}$ dominates that of $\hat{\beta}_{\rm L}$
everywhere, except for very large deviations (large $t$). It is the
heavy tails of $\hat{\beta}_{\rm M}$ that make its mean square error
infinite, otherwise it tends to be closer to $\beta$ than its rival
$\hat{\beta}_{\rm L}$.

Anderson 1976 remarks that this situation, in its extreme, resembles
the following dilemma: suppose we are estimating a parameter
$\theta$ whose true value is $\theta^\ast \approx 0$, and we have to
choose between two estimates: one, $\hat{\theta}_1$, has Cauchy
distribution, and the other, $\hat{\theta}_2$, has a normal
distribution with mean 100 and variance 1. Would anyone prefer
$\hat{\theta}_2$ only because it has finite moments?

Thus Anderson and others build a very strong case supporting the MLE
estimate $\hat{\beta}_{\rm M}$, despite its infinite moments.
Furthermore, Gleser 1983 proves that the MLE estimate
$\hat{\beta}_{\rm M}$ is the best possible in a certain formal
sense, we refer the reader to Chen and Van Ness 1994 for a detailed
survey.

\medskip \textbf{Fitting circles.}
Now we return to the circle fitting problem. It allows an
alternative approach: instead of minimizing geometric distances
\eqref{Fmain}--\eqref{argmin} one can minimize the so-called
`algebraic distances':
\beq  \label{Kasa}
   (\ha_0, \hb_0, \hR_0) = \, {\rm argmin}\, \sum_{i=1}^n \bigl[
       (x_i-a)^2+(y_i-b)^2 - R^2\bigr]^2.
\eeq
By changing parameters $A=-2a$, $B=-2b$, and $C=a^2+b^2-R^2$ one
reduces \eqref{Kasa} to a \emph{linear} least squares problem
\beq  \label{Kasa1}
   (\ha_0, \hb_0, \hR_0) = \, {\rm argmin}\, \sum_{i=1}^n \bigl[
       x_i^2+y_i^2 +Ax_i+By_i+C\bigr]^2
\eeq
which has a unique and explicit solution. This approach is known as
a \emph{simple algebraic fit} (see Chernov and Lesort 2005) or
Delogne-K{\aa}sa method (Zelniker and Clarkson 2006); it was
introduced in the 1970s. It has an obvious advantage of simplicity
over the geometric fit, which requires iterative numerical schemes.

The competition between the geometric and algebraic circle fits is
now over 30 years old, and so far it was focused on simplicity
versus accuracy. Geometric estimates $(\ha, \hb, \hR)$ are widely
known to be extremely accurate in practical applications, despite
their slight tendency to overestimate the circle's radius (the
latter was pointed out by Berman 1989). On the other hand, the
Delogne-K{\aa}sa estimates are heavily biased toward smaller
circles, see Chernov and Lesort 2004 and 2005 and references
therein, and generally much less accurate than the geometric
estimates.

Now this competition acquires a new, purely statistical momentum.
Recently Zelniker and Clarkson 2006 proved that the Delogne-K{\aa}sa
estimates $(\ha_0, \hb_0, \hR_0)$ have finite mean values whenever
$n >3$ and finite variances whenever $n>4$. Our work shows that the
geometric estimates $(\ha, \hb, \hR)$ have infinite moments.

This competition very much resembles the one described above between
the two line slope estimates: the MLE $\hat{\beta}_{\rm M}$ and the
`classical least squares' $\hat{\beta}_{\rm L}$. It would be
interesting to further compare the two circle fits along the lines
of the cited works by Anderson, Kunitomo, Gleser, and others, but
this is perhaps a research program for distant future.

\medskip \textbf{Alternative parametrizations.} One can also say that
non-existence of moments is an artifact of a poorly chosen
parametrization, and the problem is easily remedied by changing
parameters. In the case of lines, one can replace its slope $\beta$
with the angle $\theta$ the line makes with, say, the $y$-axis. Then
the line can be described as $x \cos \theta + y \sin \theta + d =0$.
Now it is easy to check that the estimates of $\theta$ and $d$ have
finite moments. These parameters are commonly used after Anderson's
work in 1976.

In the case of fitting circles, the radius $R$ can be replaced with
the curvature $\rho = 1/R$. It is easy to check that the estimate of
$\rho$ has finite moments (up to the order $2n-3$). The center
coordinates $(a,b)$ can be replaced by, say, $c=a/R$ and $d=b/R$,
which would also have finite moments. Alternatively one can replace
them with $(q,\theta)$ defined by $a = q^{-1}\cos\theta$ and $b =
q^{-1}\sin\theta$ (V.~Clakson, private communication). All these new
parameters have finite moments.

Alternatively, one can describe circles by equation
$$
      A(x^2+y^2) + Bx + Cy + D = 0
$$
subject to constraint $B^2+C^2-4AD=1$; this was proposed by Pratt
1987. Now the parameters $(A,B,C,D)$ are defined uniquely; and using
the results of Chernov and Lesort 2005 it is easy to check that they
have finite moments.

\section{Proof of Theorem 1}

It is enough to prove our theorem for the functional model. Indeed,
then in the context of the structural model the conditional
expectations of $|\ha|$, $|\hb|$, and $\hR$ for every given
realization of $\varphi_1, \ldots, \varphi_n$ will be infinite, thus
their unconditional expectations will be infinite, too.

Next we need to make a few general remarks. First, the objective
function \eqref{Fmain} may not have a minimum. For example, if the
data points are collinear, then $\inf \cF (a,b,R) =0$, but there is
no circle that would interpolate $n>2$ distinct collinear points,
hence $\cF (a,b,R) >0$ for all $a,b,R$. In that case the best fit is
achieved by a line, which can be regarded as a `degenerate circular
arc with infinite radius'.

It is proved in Chernov and Lesort 2005 that if one poses the circle
fitting problem in this `extended sense', i.e.\ as finding a circle
\emph{or a line} which minimizes the sum of squares of distances to
the given data points, then the problem always has a solution. That
is, the best fitting contour (a circle or a line) always exists. The
solution may not be unique, though, as the global minimum of the
objective function \eqref{Fmain} can be attained simultaneously on
several distinct circles, examples are given in Chernov and Lesort
2005 and Zelniker and Clarkson 2006.

In the case of multiple solutions, any one can be selected, our
theorem remains valid for any selection. If the best fit is a line,
rather than a circle (for example, if the data are collinear), then
we can set $\ha = \hb = \hR = \infty$.

This fact by itself does not prove our theorem, of course, as the
probability of such an exceptional event is zero. It shows, however,
that in nearly collinear cases the estimates $\ha, \hb, \hR$ tend to
take arbitrarily large values, and we will explore this tendency
thoroughly.

\medskip \textbf{Simple case $n=3$.} Our argument is particularly
simple if $n=3$, and this case also illustrates our main idea.

Let 3 data points be located at $(0,0)$, $(0,-1)$ and $(x,1+y)$
where $x$ and $y$ are small, say $\max \{ |x|, |y| \} \leq h =
10^{-9}$. Note that for $n=3$ the best fitting circle simply
interpolates the three given points, so by elementary geometry $\ha
=(2+3y+y^2+2x^2)/(4x)$, in particular $|\ha| \geq 1/(3|x|)$. Since
the density of $(x_3, y_3) = (x, 1+y)$ is continuous and positive,
it has a minimum value $p_0>0$ in the rectangle $|x_3| \leq h$,
$|y_3-1| \leq h$. Therefore the conditional expectation of $|\ha|$,
when the other two points are fixed, is
$$
   E\bigl(|\ha|\,\big/ B\bigr)
   \geq p_0 \int_{1-h}^{1+h} \int_{-h}^h \frac{1}{3|x|}\, dx\, dy
   = \infty,
$$
where $ B = \{(x_1, y_1)=(0,0),\ (x_2, y_2)=(0,-1)\}$.

A similar estimate holds if the points $(0,0)$ and $(0,-1)$ are
perturbed slightly, say within a little square of size $h^2$ around
their initial positions. Now the densities of $(x_1, y_1)$ and
$(x_2, y_2)$ are also positive, so a direct integration yields
$E(|\ha|) = \infty$. It is also clear that $E(\hR)=\infty$. Rotating
our construction, say by $\pi/2$, we obtain $E(|\hb|)= \infty$, too.

\medskip \textbf{General case $n>3$.}
We modify our previous construction as follows. Let $h = 10^{-9}
n^{-2}$ (here $10^{-9}$ may be replaced with any sufficiently small
constant).

We place our first point $(x_1, y_1)$ in the `lower' square
$[-h^2,h^2] \times [-1-h^2, -1+h^2]$, then $n-2$ points $(x_i,
y_i)$, $i = 2, \ldots, n-1$, in the `central' square $[-h^2,h^2]
\times [-h^2,h^2]$, and the last point $(x_n, y_n)$ in the
(horizontally extended) `upper' rectangle $[-h,h] \times [1-h^2,
1+h^2]$.

For every fixed positions of the first $n-1$ points and the fixed
$y$-coordinate $y_n$ of the last point, we will examine how the best
fitting circle changes as the $x$-coordinate $x=x_n$ of the last
point changes from $-h$ to $h$. Let $\ha(x)$ denote the first
coordinate of the circle's center (we suppress its dependence on the
other $x_i$ and $y_i$ coordinates). If the best fit is a line (and
that line is clearly almost vertical), we set $\ha= \infty$. Since
$\ha$ is large, it is more convenient to work with $\zeta(x) =
1/\ha(x)$, which is always finite and small.

Observe that all our points $(x_i, y_i)$, $1 \leq i \leq n$, are
located in the $h^2$-vicinity of three points: $(0,0)$, $(0,-1)$,
and $(x,1)$, thus the best fitting circle (or line) passes through
the $h^2$-vicinity of these three points, too. By elementary
geometry, if $x=h$, then $\ha(x) > 1/(2h)$, hence $\zeta(h) \in (0,
2h)$. Similarly, $\ha(-h)< -1/(2h)$, hence $\zeta (-h) \in (-2h,
0)$.

As $x=x_n$ changes from $-h$ to $h$, the $\zeta(x)$ function moves
from the negative interval $(-2h, 0)$ into the positive interval
$(0, 2h)$, and it stays between $-2h$ and $2h$. All we need now is
that $\zeta(x)$ behave regularly in the following sense:

\begin{lemma}[Regularity]
For any fixed values $(x_i, y_i)$, $1 \leq i \leq n-1$, and $y_n$,
as above, the function $\zeta(x)$ is differentiable and its
derivative is bounded, i.e.\ $|\zeta'(x)| \leq D$ for some constant
$D>0$. Here $D$ may depend on $n$ and $h$ but not on the fixed
coordinates $(x_i, y_i)$.
\end{lemma}

The proof of Lemma is rather technical; it is given in Appendix.
\medskip

\emph{Proof of Theorem 1}. Due to the regularity lemma, the function
$\zeta(x)$ is continuous, hence $\zeta(x_0) =0$ for some $x_0 \in
(-h, h)$. The boundedness of the derivative $\zeta'(x)$ implies that
for any $\varepsilon>0$ if $|x-x_0| < \varepsilon$, then $|\zeta(x)|
< D\varepsilon$, hence $|\ha(x)| > 1/(D\varepsilon)$. The
conditional probability of this event (when $(x_i, y_i)$ for $i =1,
\ldots, n-1$ and $y_n$ are fixed) is $\geq p_0 \varepsilon$ with
some constant $p_0>0$, due to the positivity of the density of
$(x_n,y_n)$. Therefore again, as in the $n=3$ case, the conditional
expectation of $|\ha|$ is infinite, hence so is the unconditional
expectation due to the positivity of the densities of $(x_n,y_n)$,
$i = 1, \ldots, n-1$.

Our analysis also implies $E(\hR) = \infty$. Rotating our
construction by $\pi/2$ gives $E(|\hb|)= \infty$. \qed

\medskip \textbf{Radial model.} Our theorem can be extended to
another interesting model for the circle fitting problem proposed by
Berman and Culpin 1986 and further studied by Chernov and Lesort
2004. In this model the error vector $\be_i = (\delta_i,
\varepsilon_i)$, cf.\ \eqref{errors}, satisfies $\be_i = \xi_i
\bn_i$, where $\bn_i$ is a unit normal vector to the true circle at
the true point $(x_i^\ast, y_i^\ast)$ and $\xi_i$'s are independent
normally distributed random variables with zero mean. In other
words, the noise $(\delta_i, \varepsilon_i)$ is normal but
restricted to the radial direction (perpendicular to the circle).

To extend our theorem to this model we need to assume that there are
at least three distinct true points $(x_i^\ast, y_i^\ast)$ on the
true circumference. We outline the modifications in our argument
needed to cover this new case.

Clearly it is possible that all the data points are collinear, i.e.\
there is a line $L$ such that the probability that all the data
points lie in the $h$-vicinity of $L$ is positive for any $h>0$.
Also, for at least one data point its radial direction (on which its
distribution is concentrated) must be transversal to $L$. Let that
point be $(x_n, y_n)$. Now we can repeat our construction by moving
$(x_n, y_n)$ across $L$ and keeping all the other points fixed in a
tiny vicinity of $L$. The technical analysis only requires minor
modifications in this new case, so we omit details.

\medskip \noindent\textbf{Acknowledgement}. The author is partially
supported by NSF grant DMS-0652896.

\section*{Appendix}

Here we prove our regularity Lemma 3. First we eliminate $R$ from
the picture. The objective function \eqref{Fmain} is a quadratic
polynomial in $R$, hence it has a unique global minimum in $R$ when
the other two variables $a$ and $b$ are kept fixed, and it is
attained at
\beq
      R = R(a,b) = \frac 1n \sum_{i=1}^n \sqrt{(x_i-a)^2 + (y_i-b)^2}.
         \label{Rbarr}
\eeq
This allows us to express $\cF$ as a function of $a$ and $b$ only:
\begin{align}
      \cF(a,b) &=\sum_{i=1}^n \bigl[\sqrt{(x_i-a)^2 + (y_i-b)^2} -
      R(a,b)\bigr]^2 \nonumber \\
      &=n\bigl[\bar{z} -2a\bar{x} -
      2b\bar{y}+a^2+b^2\bigr]-n[R(a,b)]^2,
         \label{Fab}
\end{align}
where for brevity we denote $z_i = x_i^2+y_i^2$; here and on we use
standard statistical `sample mean' notation $\bar{z} = \frac 1n \sum
z_i$, $\bar{x} = \frac 1n \sum x_i$, etc.

Next we switch to polar coordinates $a=\rho\cos\theta$ and
$b=\rho\sin\theta$ in which \eqref{Fab} takes form
\begin{align}
      \tfrac 1n \cF(\rho,\theta)
      &=\bar{z} -2\rho(\bar{x}\cos\theta +
      \bar{y}\sin\theta)+\rho^2\nonumber\\&\quad-
      \Bigl[\tfrac 1n\sum\sqrt{z_i-2\rho(x_i\cos\theta +
      y_i\sin\theta)+\rho^2}\Bigr]^2.
         \label{Frhotheta}
\end{align}
Note that $\cF$ in \eqref{Fab} and \eqref{Frhotheta} denotes the
same function, though expressed in different sets of variables. We
introduce more convenient notation
$$
   u_i = x_i \cos\theta + y_i\sin\theta
   \qquad\text{and}\qquad
   v_i = -x_i \sin\theta + y_i\cos\theta
$$
(observe that $u_i^2 + v_i^2 = z_i$), so that \eqref{Frhotheta}
becomes shorter:
\beq
      \tfrac 1n \cF(\rho,\theta)
      =\bar{z} -2\rho\bar{u} +\rho^2-
      \Bigl[\tfrac 1n\sum\sqrt{z_i-2\rho u_i+\rho^2}\Bigr]^2.
         \label{Frhotheta1}
\eeq
Now we introduce another variable
\begin{align} \label{wi}
   w_i &= \rho\bigl[\sqrt{z_i-2\rho u_i+\rho^2}
   -(\rho - u_i)\bigr]\nonumber\\
   &=\frac{v_i^2}{\sqrt{1-2u_i\rho^{-1}+z_i\rho^{-2}}
   +1 - u_i\rho^{-1}}.
\end{align}
From \eqref{wi} we have $\sqrt{z_i-2\rho u_i+\rho^2} = \rho - u_i +
w_i \rho^{-1}$, hence
$$
  \tfrac 1n \sum \sqrt{z_i-2\rho u_i+\rho^2} =
  \rho - \bar{u} + \bar{w}\rho^{-1}.
$$
Now \eqref{Frhotheta1} takes form
\beq
      \tfrac 1n \cF(\rho,\theta)
      =\bar{z} -\bar{u}^2 -2\bar{w}
      +2\bar{u}\bar{w}\rho^{-1}-\bar{w}^2\rho^{-2}.
         \label{Frhotheta2}
\eeq
By elementary geometry, the (averaged) objective function $\tfrac 1n
\cF$ takes all its small values (say, all values less than $h^2/10$)
on circles and lines that pass in the $h$-vicinity of the three
basic points: $(0,0)$, $(0,-1)$ and $(0,1)$. These circles and lines
have parameters restricted to the region where $\rho
> 1/(100h)$ and $|\sin \theta| < 100h$.

We replace the large parameter $\rho$ with its reciprocal $\delta =
\rho^{-1}$ and obtain
\beq
      \tfrac 1n \cF(\delta,\theta)
      =\bar{z} -\bar{u}^2 -2\bar{w}
      +2\bar{u}\bar{w}\delta-\bar{w}^2\delta^2,
         \label{Fdeltatheta}
\eeq
where
$$
   w_i = \frac{v_i^2}{\sqrt{1-2u_i\delta+z_i\delta^2}
   +1 - u_i\delta}
$$
(recall that $u_i$'s and $v_i$'s depend on $\theta$ but not on
$\rho$).

Observe that the transformation $\theta \mapsto \theta + \pi$ and
$\delta \mapsto -\delta$ leaves $w_i$'s and $\cF( \delta, \theta)$
unchanged; thus we can let $\delta$ take (small) negative values but
keep $\theta$ close to 0. More precisely, we can restrict our
analysis to the region
\beq  \label{region}
   \Omega = \bigl\{|\delta| \leq 100h \quad\text{and}\quad
   |\theta| \leq 100h\bigr\}.
\eeq
Now one can easily see that the function $\cF(\delta,\theta)$ in
$\Omega$ is regular in the following sense: it is continuous and has
bounded first and second derivatives (including partial derivatives)
with respect to its variables $\delta$ and $\theta$ and with respect
to $x = x_n$. We denote the first derivatives by $\cF_{\delta}$,
$\cF_{\theta}$, $\cF_x$ and second derivatives by
$\cF_{\delta\delta}$, $\cF_{\delta\theta}$, etc.

All these derivatives are uniformly bounded by a constant $M>0$ that
may depend on $n$ and $h$ but not on the other point coordinates.

By direct differentiation of $\cF(\delta,\theta)$ we see that
\beq \label{three}
    \cF_{\delta\delta}=
    1-\tfrac 2n+\chi_1,\quad
    \cF_{\theta\theta}=
    4-\tfrac 8n+\chi_2,\quad
    \cF_{\delta\theta}=
    \chi_3,
\eeq
where $\chi_i$ are various small quantities (that can be made as
small as we please by further decreasing $h$). Thus, $\cF$ is a
convex function that has exactly one minimum in $\Omega$ and no
other critical points.

Let $(\hat{\delta}, \hat{\theta})$ denote that unique minimum.
Differentiating equations
$$
   \cF_{\delta}(\hat{\delta}, \hat{\theta})=0
   \qquad\text{and}\qquad
   \cF_{\theta}(\hat{\delta}, \hat{\theta})=0
$$
with respect to $x$ gives
\begin{align*}
   \cF_{\delta\delta}(\hat{\delta}, \hat{\theta})\,\hat{\delta}'
   +\cF_{\delta\theta}(\hat{\delta},
   \hat{\theta})\,\hat{\theta}'
   +\cF_{\delta x}(\hat{\delta}, \hat{\theta})&=0\\
   \cF_{\theta\delta}(\hat{\delta}, \hat{\theta})\,\hat{\delta}'
   +\cF_{\theta\theta}(\hat{\delta}, \hat{\theta})\,\hat{\theta}'
   +\cF_{\theta x}(\hat{\delta}, \hat{\theta})&=0,
\end{align*}
where $\hat{\delta}'$ and $\hat{\theta}'$ denote the derivatives
with respect to $x$.

Since all partial derivatives are uniformly bounded by $M$ and the
determinant is $\approx 4-\tfrac 8n$ due to \eqref{three}, we have
that $|\hat{\delta}'| \leq 2M$ and $|\hat{\theta}'| \leq 2M$.
Lastly, recall that $\zeta = 1/\ha = \hat{\delta}/\cos
\hat{\theta}$, hence
$$
  |\zeta'| = \biggl|\frac{\hat{\theta}' \sin \hat{\theta}}
  {\cos^2 \hat{\theta}}\, \hat{\delta}
  +\frac{\hat{\delta}'}{\cos \hat{\theta}}\biggr|
  \leq 4M,
$$
which proves the lemma with $D=4M$. \qed \bigskip

\bigskip
\centerline{\sc References}
%\begin{thebibliography}{99}

%\bibitem{An81}
\medskip\noindent
Anderson D A (1981) The circular structural model.
\emph{J. R. Statist. Soc. B}, \textbf{27}, 131--141.

%\bibitem{An76}
\medskip\noindent
Anderson T W (1976) Estimation of Linear Functional Relationships:
Approximate Distributions and Connections with Simultaneous
Equations in Econometrics. \emph{J. R. Statist. Soc. B},
\textbf{38}, 1--36.

%\bibitem{An84}
\medskip\noindent
Anderson T W (1984) Estimating Linear Statistical Relationships.
\emph{Annals Statist.}, \textbf{12}, 1--45.

%\bibitem{AS82}
\medskip\noindent
Anderson T W and Sawa T (1982) Exact and Approximate Distributions
of the Maximum Likelihood Estimator of a Slope Coefficient. \emph{J.
R. Statist. Soc. B}, \textbf{44}, 52--62.

%\bibitem{BC86}
\medskip\noindent
Berman M and Culpin D (1986) The statistical behaviour of some least
squares estimators of the centre and radius of a circle. \emph{J. R.
Statist. Soc. B}, \textbf{48}, 183--196.

%\bibitem{Be89}
\medskip\noindent
Berman M (1989) Large sample bias in least squares estimators of a
circular arc center and its radius. \emph{Comp. Vision Graph. Image
Process.}, \textbf{45}, 126--128.

%\bibitem{Ch65}
\medskip\noindent
Chan N N (1965) On circular functional relationships. \emph{J. R.
Statist. Soc. B}, \textbf{27}, 45--56.

%\bibitem{CK06}
\medskip\noindent
Cheng C-L and Kukush A (2006) Non-existence of the first moment of
the adjusted least squares estimator in multivariate
errors-in-variables model, Metrika \textbf{64}, 41--46.

%\bibitem{CN94}
\medskip\noindent
Cheng C-L and Van Ness J W (1994) On Estimating Linear Relationships
when Both Variables are Subject to Errors. \emph{J. R. Statist. Soc.
B}, \textbf{56}, 167--183.

%\bibitem{CL2}
\medskip\noindent
Chernov N and Lesort C (2004) Statistical efficiency of curve
fitting algorithms. \emph{Comp. Stat. Data Anal.}, \textbf{47},
713--728.

%\bibitem{CL1}
\medskip\noindent
Chernov N and Lesort C (2005) Least squares fitting of circles.
\emph{J. Math. Imag. Vision}, \textbf{23}, 239--251.

%\bibitem{Gl83}
\medskip\noindent
Gleser L J (1983) Functional, structural and ultrastructural
errors-in-variables models. In \emph{Proc. Bus. Econ. Statist. Sect.
Am. Statist. Ass.}, pp.\ 57--66. Alexandria, VA.

%\bibitem{Ka98}
\medskip\noindent
Kanatani K (1998) Cramer-Rao lower bounds for curve fitting.
\emph{Graph. Models Image Process.}, \textbf{60}, 93--99.

%\bibitem{KMH04}
\medskip\noindent
Kukush A Markovsky I and Van Huffel S (2004) Consistent estimation
in an implicit quadratic measurement error model. Computational
Statistics and Data Analysis \textbf{47}, 123--147.

%\bibitem{Ku80}
\medskip\noindent
Kunitomo N (1980) Asymptotic Expansions of the Distributions of
Estimators in a Linear Functional Relationship and Simultaneous
Equations, \emph{J. Amer. Statist. Assoc.}, \textbf{75}, 693--700.

%\bibitem{Pr87}
\medskip\noindent
Pratt V (1987) Direct least-squares fitting of algebraic surfaces,
{\em Computer Graphics} {\bf 21}, 145--152.

%\bibitem{ZC06}
\medskip\noindent
Zelniker E E and Clarkson V L (2006) A statistical analysis of the
Delogne-K{\aa}sa method for fitting circles. \emph{Digit. Signal
Process.}, \textbf{16}, 498--522.

%\end{thebibliography}

\end{document}